\documentclass[conference]{IEEEtran}
\usepackage{color}
\usepackage{multirow}
\usepackage{cite}
\usepackage{balance}
\usepackage{amsmath}
\usepackage{amssymb}
\usepackage{enumitem}
\usepackage{arydshln}
\usepackage{fancyhdr}
\usepackage{marvosym}
\setlist[itemize]{noitemsep, topsep=0pt}
\usepackage{pdfsync}
\ifCLASSINFOpdf
\usepackage[pdftex]{graphicx}
\DeclareGraphicsExtensions{.pdf,.jpeg,.png}
\else
\usepackage[dvips]{graphicx}
\DeclareGraphicsExtensions{.eps}
\fi
\graphicspath{{figs/}}
\usepackage{epstopdf}
\usepackage{booktabs}
\usepackage{amsmath}

\renewcommand{\footnoterule}{%
  \kern -3pt
  \hrule width 0.49 \textwidth height 0.5pt
  \kern 1pt
}
\usepackage[absolute]{textpos}
\TPMargin{5pt}

\newcommand{\copyrightstatement}{
    \begin{textblock}{0.5}(0.08,0.94)   
    \noindent
         \footnotesize
          978-1-7281-0407-2/19/\$31.00~\copyright2019 IEEE 
    \end{textblock}
}

\begin{document}
\copyrightstatement
\title{Local Smart Inverter Control to Mitigate the Effects of Photovoltaic (PV) Generation Variability}

\author{\IEEEauthorblockN{Rahul Ranjan Jha and Anamika Dubey}
\IEEEauthorblockA{School of Electrical Engineering and Computer Science \\
Washington State University\\
Pullman, WA}}


%


\maketitle

\begin{abstract}
The rapidly transforming electric power distribution system due to the integration of distributed energy resources (DERs) specifically roof-top photovoltaic (PV) generation may lead to significant operational challenges. In this paper, we address the challenges related to variable power generation profile from PV resources leading to voltage fluctuations at the secondary feeder level. The objective is to develop local smart inverter control methods to reduce the voltage fluctuations at distribution buses using reactive power support. Towards this goal, an approach based on the power flow measurements of the lines connected to the PV generation buses and the local PV generation measurements is proposed to obtain the required reactive power support to mitigate voltage fluctuations. The proposed method is validated using two three-phase unbalanced test systems: IEEE-123 bus and the modified R3-12.47-2 taxonomy feeder with 329-buses. It is shown that the proposed local control approach is effective in mitigating voltage violations resulting from large changes in active power injections of the PV generators and in reducing overall voltage fluctuations due to PV generation variability. Further, it is validated that the proposed local control approach results in a better performance with regard to decreasing the voltage fluctuations compared to a standard Thevenin impedance-based method.


\end{abstract}
\vspace{0.2cm}

%
\IEEEpeerreviewmaketitle

\vspace{-0.4cm}

\section{Introduction}

The intermittent and variable nature of photovoltaic (PV) generators are known to cause  significant operational challenges for the power distribution systems. One such problem of interest is the voltage fluctuations in secondary node  resulting from variable PV generation profiles. Although, voltage fluctuations are not of major concern if nodal voltage are constrained within the prespecified voltage limits, but it may be relevant if the distribution feeder is employing with some advanced operation/control schemes. For example, it is important to mitigate the variations in nodal voltages if a volt-var control approach for conservation voltage reduction (CVR) is being employed otherwise, it may affect the realized CVR benefits. Similarly, it is crucial to minimize the voltage fluctuations if they interact with local control of feeder’s legacy voltage control devices such as voltage regulators and/or capacitor banks.  The rapidly varying voltage profile due to variable PV generation leads to unnecessary switching operations for legacy voltage control devices and thereby  loss in equipment life. This calls for effective operational scheme to manage the impacts of PV variability on power distribution systems especially on nodal voltage fluctuations. 

Note that traditional legacy voltage control devices (i.e. cap banks and voltage regulators) are not sufficient to manage PV variability related concerns. This is because, these are mechanical devices that introduce latency (30-90 sec) in responding to a requested change in set point. Since PVs are varying at much rapid interval, coordinating legacy devices to mitigate such concerns is not feasible. A new voltage control device is required that can quickly respond to PV generation variability is required to mitigate the rapid fluctuations in nodal voltages \cite{chertkov2009distributed}. To address these and related concerns, in literature the control of PV smart inverters is proposed as a viable mechanism. The modern PV systems are equipped with smart inverters that is able to provide voltage support by controlling reactive power. Smart inverters, being power electronic devices, respond almost instantaneously to a control signal and dispatch the required reactive power support.    
 
In literature, smart inverter-based control has been thoroughly employed for power distribution systems. The existing methods for smart inverter-based voltage control can be broadly divided into three categories: (i) methods based on optimal power flow (OPF) using a centralized control architecture. Here, it is assumed that there is a perfect two-way communication between the control center and smart inverters \cite{paper1,garcia2014combined,MINLPCVR}; (ii) distributed control methods that requires communication only among neighbors \cite{li2014real,kraning2014dynamic,vsulc2014optimal,robbins2012two}; and (iii) local control methods  that do not require any external communication and use only local measurements to generate control signals. Note that a local control approach is similar to primary response/control and therefore is the fastest approach and is most amenable for reducing voltage fluctuations. The existing methods for local control of smart inverters are mostly based on sensitivity matrix and volt-var droop characteristics of smart inverters \cite{shivashankar2016mitigating, Jensen, sensitivitybased, Zhao, Chamana, FarivarChe, EAbdelkarim, Liu,bokhari2016combined, zhang2018novel,zhu2016fast}. In volt-var droop control approach, a droop curve representing the mapping between the nodal voltage and reactive power dispatch is specified. The controller then autonomously responds to local voltage measurements and dispatches the reactive power as specified by the volt-var curve. Usually, the set points for the droop curve are predetermined to simply mitigate any overvoltage or undervoltage concerns \cite{shivashankar2016mitigating, sensitivitybased, Chamana}. These methods may cause voltage to oscillate, lead to a steady state error, and do not meet system-wide objectives. The optimality and the stability of the droop based voltage control methods is discussed in \cite{zhu2016fast}. In order to remove the oscillations and improve stability, several authors \cite{Zhao,Chamana,FarivarChe} has proposed methods to dynamically change the droop points based on the local measurements. Dynamically obtaining droop set-points at shorter time-intervals (1-minute or less) remain a challenging problem.

Another approach is to use sensitivity-based measures to calculate the required reactive power dispatch to mitigate the voltage fluctuation due to active power fluctuation at the current time step. The sensitivity based approach which utilizes $\frac{R}{X}$ ratio of the distribution system is used to obtain the required reactive power to reduce voltage fluctuation \cite{sensitivitybased}. However, in this work the distribution system is modeled as an equivalent balanced single-phase system. In another work, the coefficients of change in active and reactive power are determined by repeatedly solving the power flow equations \cite{Liu}. The obtained coefficients, however, are not accurate for all sets of operating conditions. In \cite{tanaka2010decentralised}, authors propose decentralized voltage control of DGs by generating voltage reference for each inverters using PI control. In \cite{bonfiglio2014optimal}, authors propose a method based on feedback linearization where system voltage is maintained by providing the required reactive power to the PVs. All of the above literature, tries to solve the problem locally without taking into account the efforts provided by other smart inverters connected in the system. Further, these methods are solved by converting the three-phase unbalanced system into single-phase equivalent. 

 The motivation of this  paper is to develop a local smart inverter control for the three-phase unbalanced system that takes the actions of other smart inverters into account without specifically communicating with those. A three-phase AC linear power flow based formulation is proposed and reactive power dispatch is obtained in every 1-min interval by each smart inverter using local measurements for change in active power injected by the PV at a node and the change in power flowing into the children nodes.

\section{Linear Three-Phase Power Flow Model}
In this section, a linear three-phase unbalanced power flow model is 
introduced. The linear three-phase power flow is based on the branch flow model where the variables are the power flow in the lines and the node voltages. 

\begin{figure}[t]
\vspace{-0.6cm}
\centering
\includegraphics[width=0.45 \textwidth]{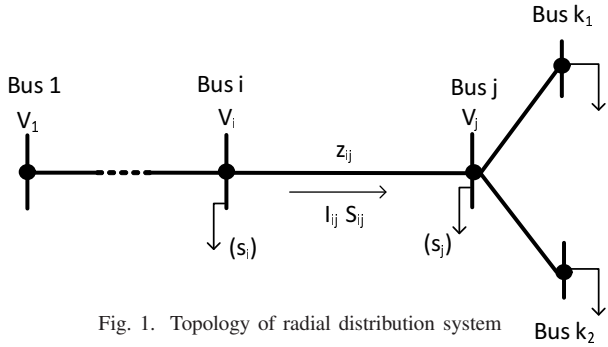}
\vspace{-1 cm}
	\caption{Topology of radial distribution system }
	\vspace{-0.9 cm}
	\label{fig:1}
\end{figure}

As the  distribution system mostly operates radially, it can be considered a directed graph. The radial distribution feeder in Fig.1 is considered as a  directed graph $\mathcal{G}(\mathcal{N, E})$ , where, $\mathcal{N}$ is the number of nodes and  $\mathcal{E}$ is the number of edges in the graph. The edge ($i,j$) connects nodes $i$ and $j$ where node $i$ is the parent of  node $j$. The complex number $z_{ij} = r_{ij} + \iota x_{ij}$ represents the complex impedance of the edge $(i,j) \in \mathcal{E}$ connecting nodes $i$ and $j$. Also, for each edge ($i,j$), assume that the apparent power flow is $S_{ij} = P_{ij} + \iota Q_{ij}$ and complex line current is $I_{ij}$. For each node $(i) \in \mathcal{N}$, let $V_i$ be the complex voltage and $s_i = p_i + \iota q_i$ be the net apparent power injection (generation minus demand). 

The  voltage drop and power balance equation for a three phase distribution system proposed in \cite{gan2014convex} are shown in equation (1)-(2).
\begin{eqnarray}
v_i  = v_j + (S_{ij}z_{ij}^H+z_{ij}S_{ij}^H) - z_{ij}I_{ij}I_{ij}^H z_{ij}^H \\
s_{L,j} = \text{diag}(S_{ij} - z_{ij}I_{ij}I_{ij}^H) -  \sum_{k:j \rightarrow k}{\text{diag}(S_{jk})}   
\end{eqnarray} 

The equations (1) and (2) are the nonlinear sets of equations due to $I_{ij}I_{ij}^H$. These equations are converted into linear sets of equations by assuming that the loss occurring in each individual edges are much smaller than the power flow in the edges. Further, to reduce the number of variables due to mutual coupling in the system. It is assumed that the phase angle difference among the phases are $\frac{2\pi}{3}$ apart \cite{gan2014convex}.

The linear power flow equations obtained for a three-phase unbalanced system is given as:
\begin{eqnarray} \label{eq:3}
  p^p_{j}   &=& P^{pp}_{ij} -\sum_{k:j \rightarrow k} P^{pp}_{jk}   \hspace{1 cm}  p \in {a,b,c}\\
  q^p_{j}   &=&Q^{pp}_{ij} - \sum_{k:j \rightarrow k} Q^{pp}_{jk} +  \hspace{1 cm}  p \in {a,b,c}\\
   v_i^p   &=& v_j^p + \sum_{q \in \phi_ j} 2\Re [S_{ij}^{pq}(z_{ij}^{pq})^*] \hspace{.4 cm}  p,q \in {a,b,c} 
\end{eqnarray}

where, $p^p_{Lj}$ and $q^p_{Lj}$ are the active and reactive power demand of the load. $P^{pp}_{ij}$ and $Q^{pp}_{ij}$ are the active and reactive power flow in the edges and $v_i^p = (V_i^p)^2$ is the square of voltage magnitude at a node $i$.

\section{Local Smart Inverter Control}
The local control is required to reduce the voltage fluctuation caused due to variability in the power injection by the PV. In this section, the local control based on equivalent Thevenin's impedance is introduced and another local control is proposed based on measurement of the power flow in the lines. The  local control methods operate at every 1-min interval and are designed to reduce the voltage fluctuation caused due to variability in the power injection by the PV. The local control provides the required reactive power from the smart inverters'. The control decisions for the smart inverters' are obtained by only providing the local measurements at the point of connection of the PV.

\subsection{Local Smart Inverter Control - using Equivalent Thevenin Impedance Method}

\begin{figure}[t]
\centering
\includegraphics[width=3.5 in]{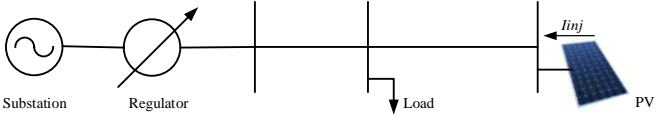}
 \vspace{-0.2cm}
	\caption{\footnotesize Distribution system connected to PV}
\vspace{-0.6cm}
	\label{fig:2}
\end{figure}

The objective of the local control is to reduce the voltage fluctuations at the nodes of the distribution system caused due to change in active power injected at the PV generation nodes.  Here, PV system is assumed as a current source injecting $I_{inj}$ at a node as shown in Fig.2. The voltage at a node $i$ for a phase $p$, $V^p_{i}$, is given by (\ref{eq6}) where, $I^p_{inj}$ is the current injected by the PV system and $Z_{i}^{p}$ is the equivalent Thevenin's impedance at node $i$ for phase $p$.  
\begin{equation}\label{eq6}
    V^p_{i} = Z_{i}^{p}I^p_{inj} = (R_{i}^{p}+\iota X_{i}^{p})I^p_{inj} \hspace{0.5 cm} p\in \{a,b,c\}
\end{equation}

The current injected by PV at a node $i$, $I^p_{inj}$, is obtained using (\ref{eq7}), where $P^p_{inj}$ and $Q^p_{inj}$ are the per-phase active and reactive power injected at node $i$ at nominal voltage level, $V^p_{nom}$. 
\begin{equation}\label{eq7}
    I^p_{inj} = \frac{(P^p_{inj}+\iota Q^p_{inj})^*}{V^p_{nom}} \hspace{0.5 cm} p\in \{a,b,c\}
\end{equation}

The voltage fluctuation at a node is dependent on the change in active and reactive power at the node. The relationship between  change in voltage, active power and reactive power is given as:
\begin{equation}
dV^p_i = \frac{\partial V^p_i}{\partial P^p_{inj}}dP^p_{inj} + \frac{\partial V^p_i}{\partial Q^p_{inj}}dQ^p_{inj}  \hspace{0.4 cm} p\in \{a,b,c\}
\end{equation}

Further, the relationship between change in voltage wrt to Thevenin's impedance is obtained as:
\begin{equation} \label{eq9}
\begin{split}
    dV^p_i = \left(\frac{R_{i}}{V^p_{nom}}dP^p_{inj}+ \frac{X_{i}}{V^p_{nom}}dQ^p_{inj}\right) \\
    + \iota \left(\frac{X_{i}}{V^p_{nom}}dP^p_{inj}- \frac{R_{i}}{V^p_{nom}}dQ^p_{inj}\right)
\end{split}
\end{equation}

The change in voltage is approximately equal to real part of the product of current and impedance \cite{kers}. Hence, by ignoring the imaginary part and putting $dV^p_i$  to zero in equation (9). The change in reactive power required to reduce the voltage fluctuation is proportional to change in active power and $R/X$ ratio of the Thevenin's impedance for the node as given in (10).

\begin{equation} \label{eq10}
    dq^p_{inj}(t) = -\frac{R^p_{i}}{X_{i}^p}dP^p_{inj}(t)
\end{equation}

The Thevenin's impedance at a node $i$ is obtained by ignoring the shunt impedance of the distribution system. Also, in this work, the impedance value obtained is assumed to be constant irrespective of change in operating condition of the system. Hence, the required reactive power obtained using this method is not accurate to reduce to the voltage fluctuation. Further, in this method the reactive power support provided by other smart inverters is ignored, means there is no coordination among the smart inverters. Thus, the obtained method is only the approximate method. Another local control method is proposed in the following section, where the effect of reactive power support from other smart inverters are taken into consideration. 

\subsection{Local Smart Inverter Control - using Power Flow Measurements}
The power flow measurement based local control is proposed to reduce the voltage fluctuation in the distribution system. The proposed local control is based on the linear three-phase power flow introduced in section II. As the proposed method includes the effects of local control actions from other smart inverters it reduces the voltage fluctuation better than the Thevenin's impedance based method. 

The change of nodal voltages due to change in PV power injection at a node can be obtained by differentiating the voltage equation (5) as described in equation (11).   
\begin{equation} \label{eq11}
\begin{split}
    dv_i^p - dv_j^p =  2 \left(r_{ij}^{p}dP_{ij}^{p} + x_{ij}^{p}dQ_{ij}^{p}\right) 
    \end{split}
\end{equation}

Also, the change in active and reactive power flow in the line is given as:
\vspace{-0.1cm}
\begin{eqnarray} 
dp^p_{j}   &=& dP^{p}_{ij} -\sum_{k:j \rightarrow k} dP^{p}_{jk}   \hspace{1 cm}  p \in {a,b,c}\\
dq^p_{j}   &=& dQ^{p}_{ij} - \sum_{k:j \rightarrow k}d Q^{p}_{jk} +  \hspace{1 cm}  p \in {a,b,c}
\end{eqnarray}

The change in  voltage as a function of change in power flow in the children nodes is obtained from equation (11), (12) and (13) as shown in equation (14).
\vspace{-0.1cm}
\begin{equation} \label{eq14}
dv_i^p - dv_j^p =  2(r_{ij}^{p}(\sum_{k:j \rightarrow k} dP^{p}_{jk}- dp^p_{j} )+ x_{ij}^{p}(\sum_{k:j \rightarrow k}d Q^{p}_{jk} - dq^p_{j}))
\end{equation}

Next, we equate the changes in nodal voltages to zero i.e., $dv_i^p = dv_j^p = 0$. We obtain the mathematical relation between the change in reactive power flow as a function of change in active power flow for each phase as:
\vspace{-0.1cm}
\begin{equation} \label{eq15}
\begin{split}
dq^p_{j} =   \frac{r_{ij}^{p}}{x_{ij}^{p}}\left ((\sum_{k:j \rightarrow k} dP^{p}_{jk}- dp^p_{j} )+ \sum_{k:j \rightarrow k}d Q^{p}_{jk}\right) p \in {a,b,c}
\end{split}
\end{equation}

It is to be noted that the inverters  installed downstream from a PV node are providing reactive power in order to reduce the voltage fluctuations. These reactive power support will change the power flow measurements at the downstream edges. Therefore, in equation (15)  the power flow measurements at the downstream nodes are  included to realize the reactive power support provided by other inverters in the system. Thus, the proposed local control method to reduce the nodal voltage fluctuation in the distribution system is better compared to the Thevenin impedance-based method. The implementation of the proposed local control in the smart inverters requires only the measurements of the change in active power injected by the PV system at a node and the change in active and reactive power flowing through the children nodes.

\begin{figure}[t]
\centering
\includegraphics[width=3.8in]{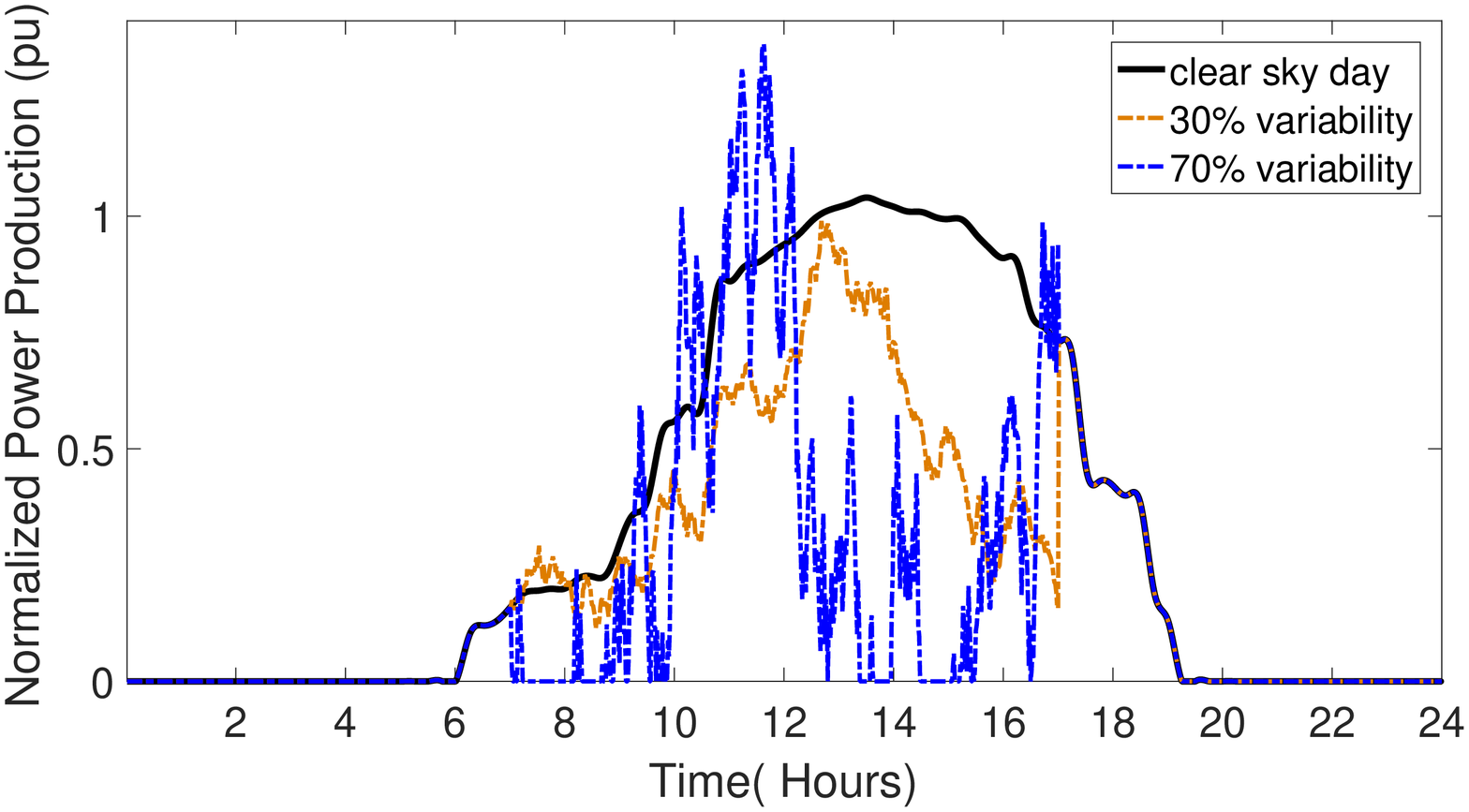}
 \vspace{-0.8cm}
	\caption{\footnotesize PV generation profile with variability.}
\vspace{-0.2 cm}
	\label{fig:3}
\end{figure}

\begin{figure}[t]
\centering
\includegraphics[width=3 in]{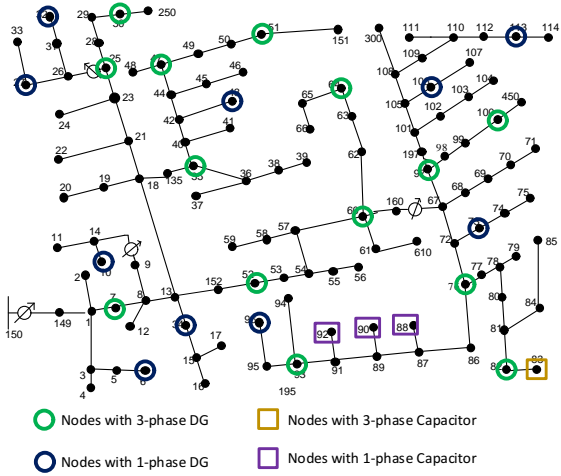}
\vspace{0.1cm}
	\caption{\footnotesize IEEE 123-bus distribution test feeder.}
\vspace{-0.5cm}
	\label{fig:4}
\end{figure}

\section{Results and Discussion}
Solar PV produces power from the irradiance following on the surface of the PV panel. Typically, for a clear sky day, it is assumed that the power generation has smooth variation and follow parabolic profile as shown in Fig.3. However, for a cloudy day, there is rapid variation in the irriadance following on the PV panel which causes fluctuation in the PV power production (see Fig.3).  To simulate variability in the PV power production, we have used the three-sigma rule  \cite{threesigma}. The The normalized PV power production with 30\%  and 70\% variability is shown in Fig. 3. 
The local control discussed in section III is validated using the modified IEEE-123 bus system \cite{TestFeeder} (see Fig.4). The feeder is modified to include a total of 55 nodes populated with PVs spread across the feeder; three-phase PVs of 172.5 kVA or 69 kVA ratings and single phase PVs of 23 kVA or 11.5 kVA ratings. Similarly,  R3-12.47-2  feeder \cite{feeder329} (see Fig.5) is modified to include a total of 50 nodes populated with PVs spread across the feeder. The three-phase PVs of ratings 690 kVA, 345 kVA, 69 kVA, 172.5 kVA, 138 kVA or 34.5 kVA and single phase PVs of ratings 23 kVA or 11.5 kVA are installed at various location in the system. The simulations are conducted using Matlab and OpenDSS and the interface between Matlab and OpenDSS is established using COM interface. The local control algorithm is executed in Matlab and the obtained control variables are provided to the distribution model simulated in OpenDSS. Therefore, although control decisions are based on linearized power flow model, OpenDSS simulates actual nonlinear three-phase system and provides a suitable environment for evaluating the performance of the proposed methods on real-world distribution systems.  All the simulations are done on i7 3.41 GHz processor with 16 GB of RAM. 
 
\vspace{-0.2 cm}
\begin{figure*}[t]
\centering
\includegraphics[width=6.3 in ]{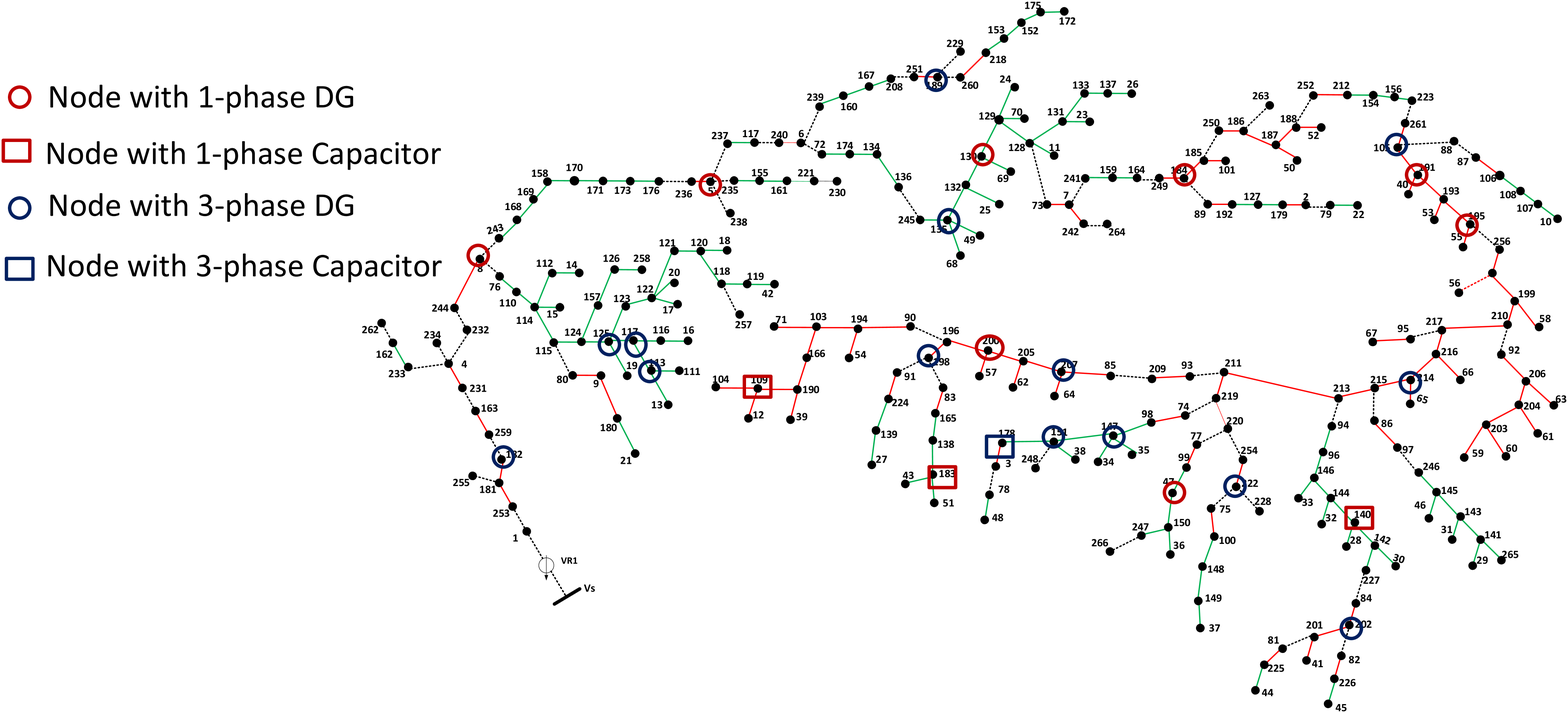}
	\caption{Modified  R3-12.47-2 distribution test feeder.}
	\label{fig:5}
	\vspace{-0.3cm}
\end{figure*}

\subsection{Reduction in Voltage Fluctuation}
In the following section, the effect of the proposed local control in reducing the voltage fluctuation is discussed. The effect of local control on voltage fluctuations is quantified by a power quality index called system average voltage fluctuation index (SAVFI) proposed in \cite{Ding}. The voltage fluctuation $ \Delta V_i(t)$ is defined as the average of the difference in the voltage magnitude $ V_i$ between the two consecutive time interval.

\begin{equation}
    \Delta V_i = |V_i^{t+1}-V_i^t|
\end{equation}

The SAVFI is defined as the average of voltage fluctuation for a time interval $T$ (here T = 15 min). 

\begin{equation}
  SAVFI = \frac{1}{T} \sum_{n=1}^{T} \Delta V_i(t)
\end{equation}

The SAVFI value for selected nodes of the IEEE 123-bus system is shown in Table I for 30\% and 70\% PV output power variability cases. It can be observed from the table the SAVFI increases as PV variability is increased. Also, it is observed that as the distance of the nodes from the substation increases the SAVFI value increases for same levels of PV variability. In order to reduce the voltage fluctuation within the 15-min time interval, both the impedance-based and the proposed power flow measurement-based local control methods are implemented at each PV location. It can be observed from the SAVFI value (see Table I), that both  the local control are able to reduce the voltage fluctuations at a node. Also, the proposed power flow measurement-based local control method is able to reduce the voltage fluctuations more efficiently as compared to impedance-based method. This is because, the proposed power flow measurement based method takes  the reactive power support provided by other smart inverters into account. 

\begin{table}[t]
		\centering
	     \caption{SAVFI for the IEEE 123-bus system at different variability level}
	     	\vspace{-0.2cm}
		\label{singletable}
		\begin{tabular}{c|c|c|c}
		\hline
		\hline
			Node  & Without local & Impedance-based & Proposed local  \\
			number & control & control & control  \\
            \hline
            \multicolumn{4}{c}{30 \% PV power variability}\\
            \hline
			7 & 2.019 & 1.12 & 0.839 \\
			\hline
			35 & 3.973 & 1.2 & 0.842 \\
			\hline
			64 & 5.505 & 1.63 & 1.116 \\
			\hline
			83 & 5.714 & 1.4 & 0.997 \\
            \hline
            114 & 5.87 & 1.55 &1.18 \\
            \hline
            \multicolumn{4}{c}{70 \% PV power variability}\\
            \hline
            7 & 2.94 & 1.76 & 1.34 \\
			\hline
			35 & 5.83 & 2.27 & 1.63 \\
			\hline
			64 & 7.87 & 2.42 &2.08 \\
			\hline
			83 & 8.15 & 2.95 & 1.62 \\
            \hline
            114 & 8.38 & 3.13 & 2.38 \\
	\hline
	\hline
\end{tabular}
\vspace{-1.3cm}
\end{table}

Similarly, the proposed local control methods are demonstrated to reduce voltage fluctuation in the  modified  R3-12.47-2 in Table II. From Table II, it can be observed that the SAVFI is higher for the nodes away from the substation. It can also be verified from the table that the proposed local control based on power flow measurements is relatively more effective in reducing the voltage fluctuations for different levels of PV generation variability. 

\begin{table}[t]
		\centering
		\caption{SAVFI for the  modified  R3-12.47-2 system at different variability level}
		\label{singletable1}
		\begin{tabular}{c|c|c|c}
		\hline
		\hline
		Node  & Without local & Impedance-based & Proposed \\
		number	& control & control & control \\
            \hline
            \multicolumn{4}{c}{30 \% PV output power variability}\\
            \hline
			5 & 2.844 & 0.77 & 0.577 \\
			\hline
			98 & 11.67 & 4.01 & 3.72 \\
			\hline
			105 & 11.32 & 4.13  & 3.77 \\
            \hline
            151 & 11.86 & 4.2 & 3.84 \\
            \hline
            214 & 11.65 & 4.31 & 3.93 \\
            \hline
            \multicolumn{4}{c}{70 \% PV output power variability}\\
            \hline
            5 & 3.33 & 0.79 &0.59  \\
            \hline
            98 & 12.12 & 4.1 & 3.85 \\
            \hline
            105 & 11.79 & 4.13 & 3.78\\
            \hline
            151 & 12.14 &4.22 & 3.86 \\
            \hline
            214 & 11.92 & 4.33 & 3.95 \\
	\hline
	\hline
\end{tabular}
\vspace{-0.8 cm}
\end{table}

\subsection{Mitigating Voltage Violations}
The over-voltage condition at nodes in the distribution system will be observed when the load demand is less and the actual PV power production is increased above the forcasted PV power production. The Fig.6, shows the over-voltage scenario in the IEEE-123 node system at the node 7. It can be observed from the figure that the voltage at nodes 7  are above 1.05 pu for few time steps.  This violates the specified ANSI voltage limits (0.95-1.05). The local control methods are effective in reducing the voltage fluctuation as well as mitigating the voltage violation at node 7.

Similarly, the under-voltage condition at a node is observed when the load demand is high and the predicted PV power production is larger than the actual power production. In Fig.7 and Fig. 8, the under-voltage scenarios are simulated for the IEEE 123-bus and the  modified  R3-12.47-2 test feeders, respectively.  For the IEEE-123 node system, node 114 observe voltages below 0.95 pu (see Fig.7 ) for few time step. The proposed local control is able to eliminate the under-voltage conditions at node 114. Also, it can observed that the power flow measurement based local control is able to reduce the voltage fluctuation more effectively than the impedance based method. Similarly, the under-voltage scenario is created for the  modified  R3-12.47-2 system and it can observed from Fig. 8 that nodes 222 observe voltages less than 0.95 pu. The proposed local control is able to mitigate the under-voltage scenario as well as able to reduce the voltage fluctuation at the node. 

\begin{figure}[t]
\centering
\includegraphics[width=0.5\textwidth]{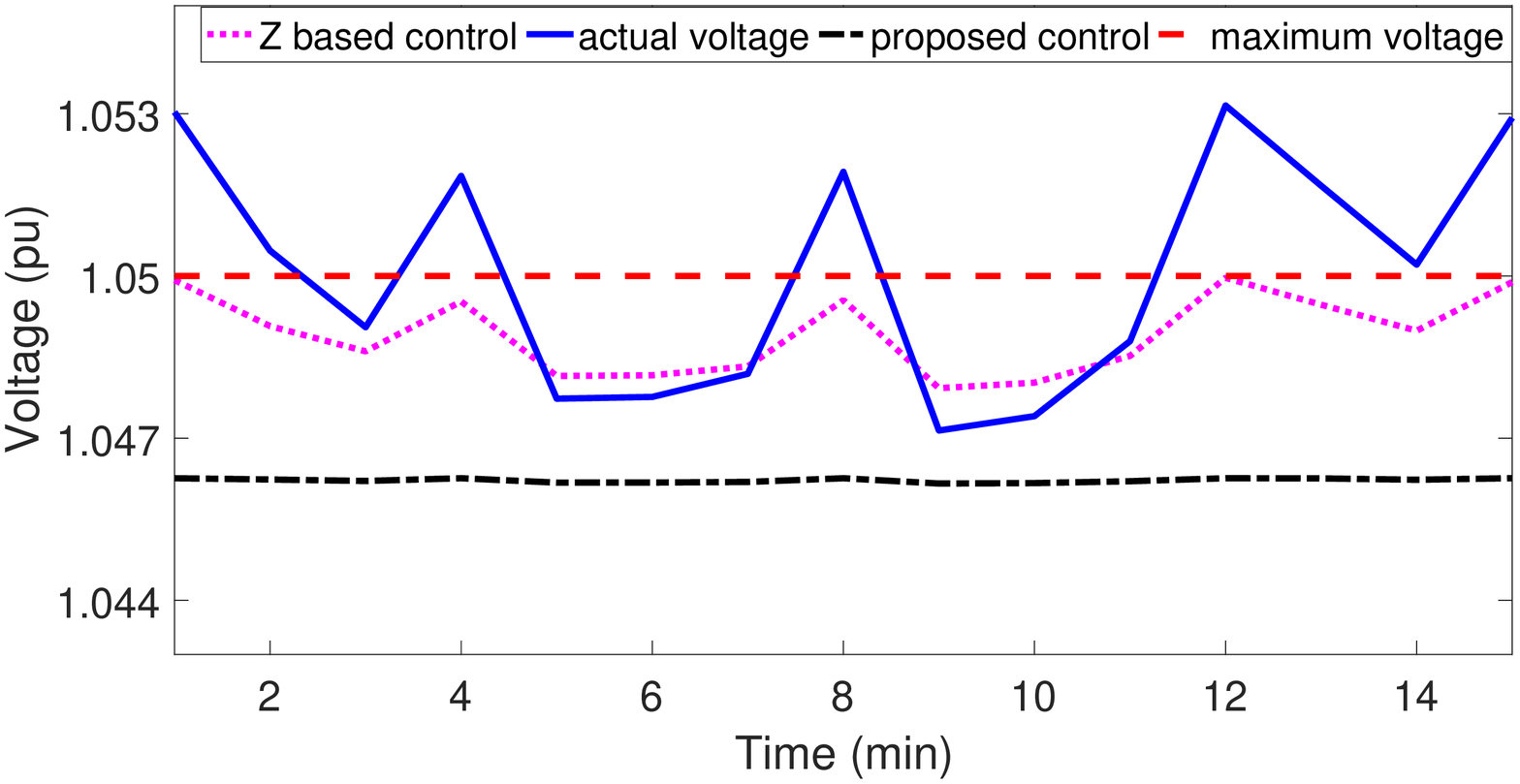}
\vspace{-0.6 cm}
	\caption{\footnotesize  Over-voltage  scenario for the IEEE 123-bus system}
\vspace{-1.10 cm}
	\label{fig:6}
\end{figure}

\vspace{-0.3cm}
\section{Conclusions}
The variable power generation resources such as roof-top PVs may lead to nodal voltage fluctuations. In this paper, a local control approach is proposed using smart inverter control that uses local measurements to calculate the required reactive power dispatch for a measured change in active power generation. The objective of the proposed local control approach is  to reduce the voltage fluctuations as well as to mitigate the over-voltage and under-voltage conditions for nodal voltages resulting from generation variability. The local control approach is based on local PV generation measurements and the power flow measurements to the children nodes. It is demonstrated that the proposed local control is effective in reducing the voltage fluctuations and also able to mitigating any voltage violation concerns. Finally, the proposed power flow measurement-based local control approach is shown to be relatively more effective compared to the Thevenin's equivalent impedance-based method. 
 
\vspace{-0.1 cm}
 \begin{figure}[t]
\centering
\includegraphics[width=0.51\textwidth]{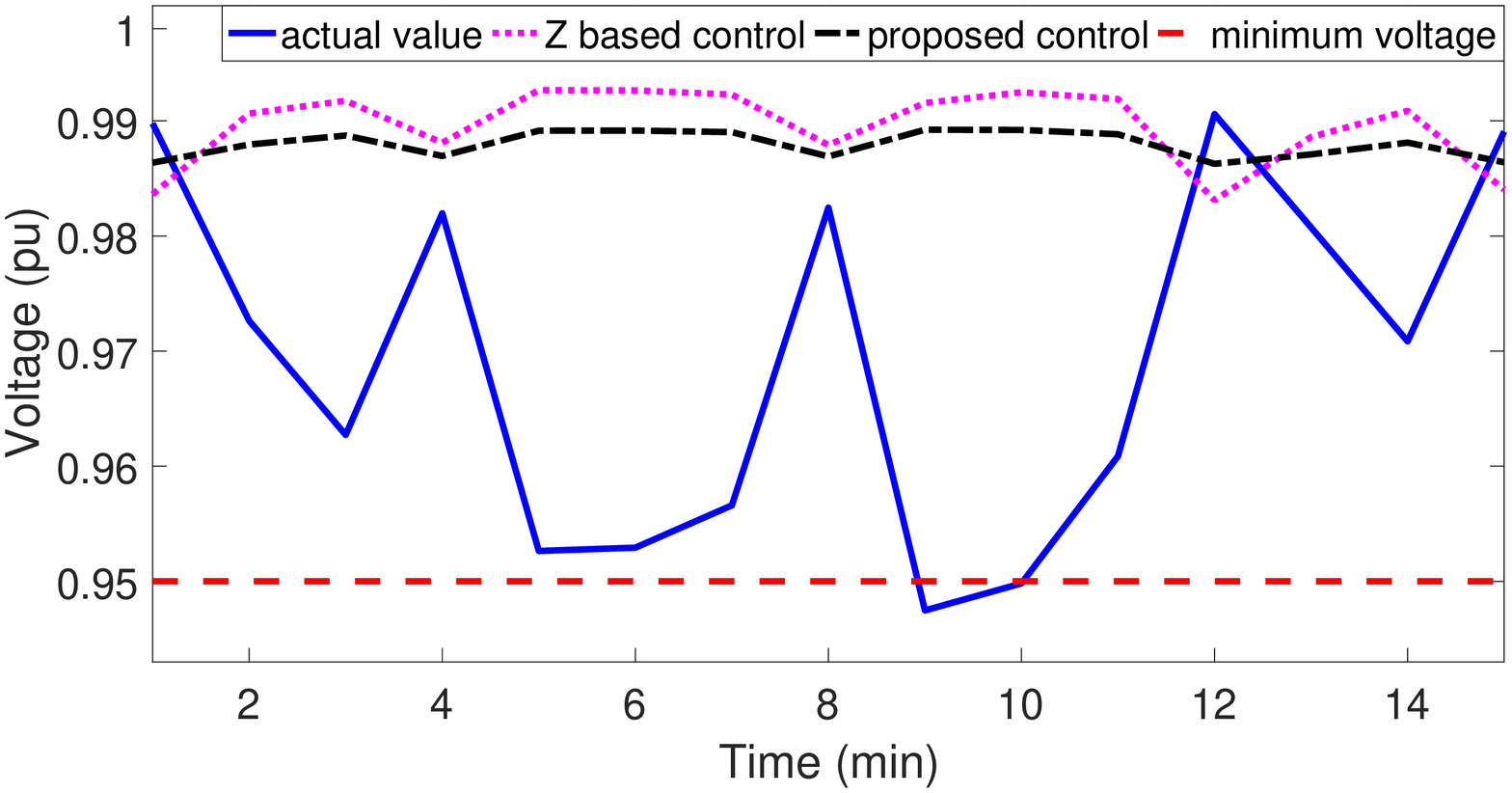}
\vspace{-0.8 cm}
	\caption{\footnotesize  Undervoltage  scenario for the IEEE 123-bus system}
\vspace{-0.5 cm}
	\label{fig:7}
\end{figure}

\begin{figure}[t]
\centering
\includegraphics[width=0.51\textwidth]{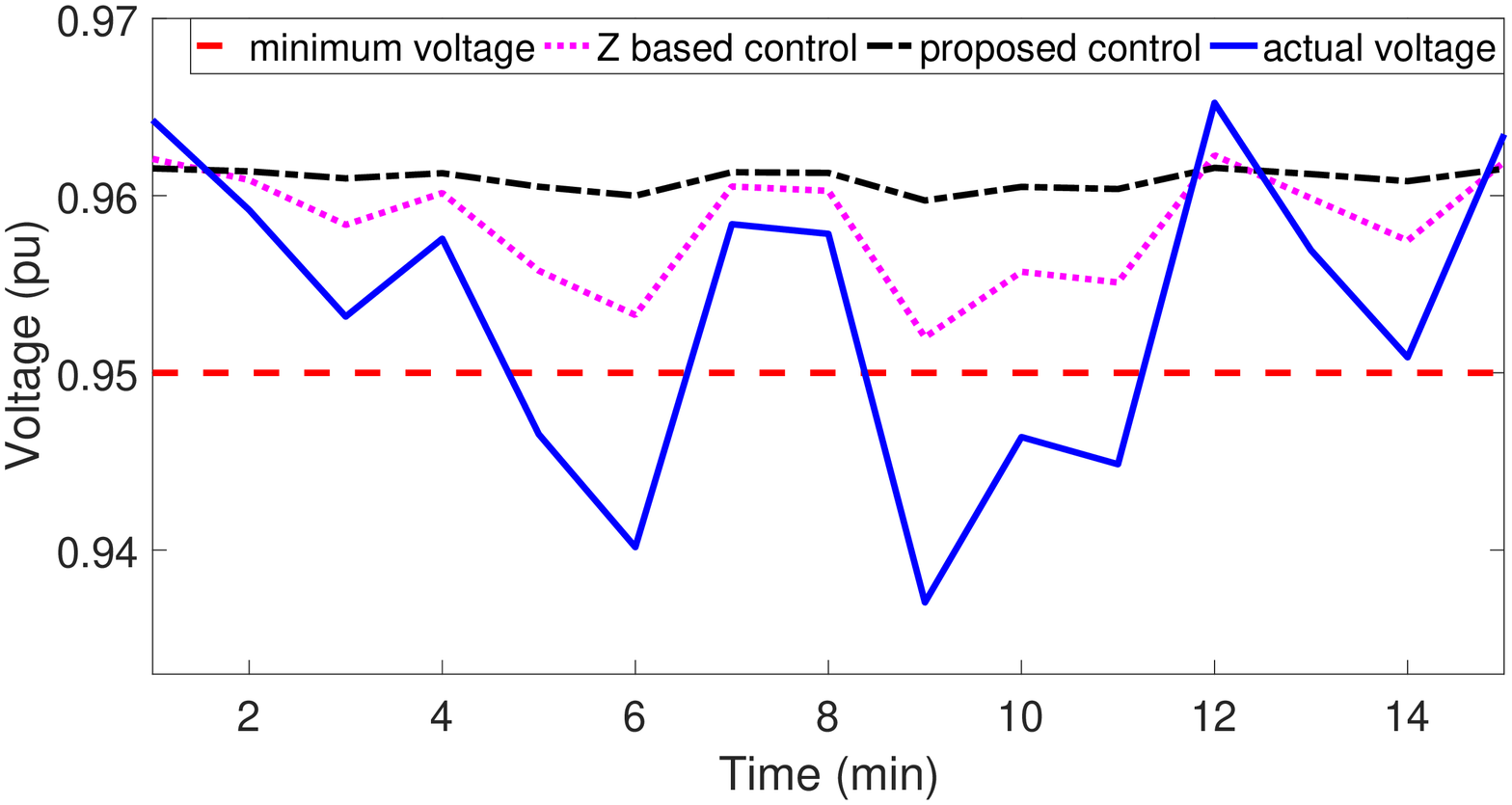}
\vspace{-0.8 cm}
	\caption{\footnotesize  Undervoltage  scenario for the modified  R3-12.47-2 system}
\vspace{-1.6 cm}
	\label{fig:8}
\end{figure}



\vspace{-0.1 cm}
\bibliographystyle{IEEEtran}
%
%
  
\bibliography{references}


\end{document}